\documentclass[leqno, 10pt, draft]{article}
\usepackage[cp1251]{inputenc}
\usepackage{ amsfonts, amssymb, amsmath, amsthm}
\usepackage[ukrainian, russian, english, french]{babel}

\renewcommand{\normalsize}{\fontsize{14}{20}\selectfont}
\begin{document}
\selectlanguage{english} \thispagestyle{empty}
 \pagestyle{myheadings}              



\vskip 2mm

 \begin{center}
 {\LARGE \bf  Approximation by Fourier sums on the classes of generalized Poisson integrals }
 \end{center}

 \begin{center}
\noindent {\rm A.\,S. SERDYUK, T.\,A. STEPANYUK}
 \end{center}

 \begin{center}
  Institute of Mathematics of Ukrainian National Academy of Sciences, 3, Tereshchenkivska st., 01601, Kyiv-4, Ukraine \\ {serdyuk@imath.kiev.ua, tstepaniuk@imath.kiev.ua}
 \end{center}

\vskip 3.5mm

{\rm  \large A\,b\,s\,t\,r\,a\,c\,t  \ We present a survey of results related to the solution of Kolmogorov--Nikolsky problem for Fourier sums  on the classes of generalized Poisson integrals
  $C^{\alpha,r}_{\beta,p}$,  which consists in finding of asymptotic equalities for exact upper boundaries of uniform norms of deviations of partial Fourier sums on the classes of  $2\pi$--periodic functions  $C^{\alpha,r}_{\beta,p}$, which are defined as convolutions of the functions, which belong to the unit balls pf the spaces $L_{p}$, $1\leq p\leq \infty$, with generalized Poisson kernels 
 $$
P_{\alpha,r,\beta}(t)=\sum\limits_{k=1}^{\infty}e^{-\alpha k^{r}}\cos
\big(kt-\frac{\beta\pi}{2}\big), \ \alpha>0, r>0, \  \beta\in
    \mathbb{R}.
$$
\hfill}

\thispagestyle{empty} \normalsize \vskip 3.5mm


In present paper, we give a survey of current results about best possible estimates of norm of deviations of partial Fourier sums on the sets of generalized Poisson integrals in the uniform and integral metric.

The partial Fourier sums
\begin{equation*}
S_{n}(f)=S_{n}(f;x)=
\frac{a_{0}}{2} +\sum\limits_{k=1}^{n} a_{k}\cos kx+b_{k}\sin kx,
\end{equation*}
\begin{equation*}
a_{k}=a_{k}(f)=\frac{1}{\pi} \int\limits_{-\pi}^{\pi}f(t)\cos kt dt, \ \ k=0,1,2,...,
\end{equation*}
\begin{equation*}
b_{k}=b_{k}(f)=\frac{1}{\pi} \int\limits_{-\pi}^{\pi}f(t)\sin kt dt, \ \ k=1,2,...
\end{equation*}
 occupy a prominent place among linear methods of polynomial trigonometric approximation, due to their simplicity and good approximative properties.
 Since the beginning of the 20th century, the problem of investigating the speed of convergence of the Fourier series, depending on the structural and smoothness properties of the function, has been studied by many mathematicians.

Results, related to certain approximate properties of Fourier sums can be found in fundamental monographs \cite{Akhiezer}--\cite{Dzyadyk}, \cite{Korn, Step monog 1987, Stepanets1,  Timan_AF1960, Zigmund} and other.

Let $L_{p}$,
$1\leq p<\infty$, be the space of $2\pi$--periodic functions $f$ summable to the power $p$
on  $[0,2\pi)$, in which
the norm is given by the formula
$\|f\|_{p}=\Big(\int\limits_{0}^{2\pi}|f(t)|^{p}dt\Big)^{\frac{1}{p}}$; $L_{\infty}$ be the space of measurable and essentially bounded   $2\pi$--periodic functions  $f$ with the norm
$\|f\|_{\infty}=\mathop{\rm{ess}\sup}\limits_{t}|f(t)|$; $C$ be the space of continuous $2\pi$--periodic functions  $f$, in which the norm is specified by the equality
 ${\|f\|_{C}=\max\limits_{t}|f(t)|}$.
 
 If we consider the estimates of the deviation of norms of the Fourier sums in the uniform metric, then the Lebesgue inequality is the  fundamental result regarding the estimate of the approximate properties of the mentioned sums  (see, e.g., \cite{Lebesque, Natanson, Step monog 1987, Stepanets1})

 \begin{equation}\label{SerdyukStepanyuk_FilomatLebeqIneq}
\| f(\cdot)-S_{n-1}(f;\cdot) \|_{C} \leq (1+ L_{n-1})E_{n}(f)_{C}, \ f\in C, \  \ n\in\mathbb{N},
\end{equation}
where
 \begin{equation*}
L_{n-1}=\frac{1}{\pi}\int\limits_{-\pi}^{\pi}|D_{n-1}(t)|dt,
\end{equation*}
  are Lebesgue constants of Fourier sums,
\begin{equation*}
D_{n-1}(t):=
\frac{1}{2}+\sum\limits_{k=1}^{\infty}\cos kt=\frac{\sin(n-\frac{1}{2})t}{2\sin \frac{t}{2}}
\end{equation*}
are Dirichlet kernels, and
 \begin{equation*}
E_{n}(f)_{C}= \inf\limits_{T_{n-1}\in \mathfrak{T}_{2n-1}} \| f- T_{n-1} \|_{C}
\end{equation*}
 is the best uniform approximation of the function  $f$ by trigonometric polynomials of the order $n-1$.
 
 For the Lebesgue constants  $L_{n}$ the following asymptotic equality takes place as  $n\rightarrow\infty$ (see \cite{Fejer})
 \begin{equation}\label{LebesgueConstant}
L_{n}=\frac{4}{\pi^{2}} \ln n+ \mathcal{O}(1), \ \ n\rightarrow\infty.
\end{equation} 

More exact estimates for the differences  $L_{n}-\frac{4}{\pi^{2}} \ln (n+a)$, where $a>0$, can be found in the works  \cite{Akhiezer, Dzyadyk, Galkin, Natanson, shakirov2018, ZhukNatanson} and others. 

Taking into account \eqref{LebesgueConstant}, the inequality (\ref{SerdyukStepanyuk_FilomatLebeqIneq}) can be written in the form 
\begin{equation}\label{SerdyukStepanyuk_FilomatLebeqIneq0}
\| f(\cdot)-S_{n-1}(f;\cdot) \|_{C} \leq \left(\frac{4}{\pi^{2} }\ln n +\mathcal{O}(1) \right)E_{n}(f)_{C}.
\end{equation}

On the whole space  $C$ the inequality (\ref{SerdyukStepanyuk_FilomatLebeqIneq0}) is asymptotically exact. At the same there exist subsets of functions from  $C$ and for elements of these subsets the inequality  (\ref{SerdyukStepanyuk_FilomatLebeqIneq0}) is not exact even by order.



Existence of such subsets follows, in particular, from Oskolkov results  \cite{Oskolkov}, who proved that for any $f\in C$  the following inequality holds
 \begin{equation}\label{oskolkovEstimate}
\| f(\cdot)-S_{n-1}(f;\cdot) \|_{C} \leq K \sum\limits_{v=n}^{2n-1}\frac{E_{v}(f)_{C}}{v-n+1}, \ f\in C,
\end{equation}
where $K$ is some absolute constant and he showed that on the classes  $$C(\varepsilon):=\{f\in C: \ E_{v}(f)_{C}\leq\varepsilon_{v}, \ v\in\mathbb{Z}_{+} \},
$$
generated by sequences  $\varepsilon=\{ \varepsilon_{k}\}_{k=0}$ of nonnegative numbersm such that  $\varepsilon_{k}\downarrow 0$, $k\rightarrow\infty$, the estimate \eqref{oskolkovEstimate} is exact by the order, i.e.
\begin{equation}\label{oskolkovEstimate_1}
K_{1}\sum\limits_{v=n}^{2n-1}\frac{\varepsilon_{v}}{v-n+1},\leq \sup\limits_{f \in C(\varepsilon)}\| f(\cdot)-S_{n-1}(f;\cdot) \|_{C} \leq K_{2} \sum\limits_{v=n}^{2n-1}\frac{\varepsilon_{v}}{v-n+1}, 
\end{equation}
where $K_{1}$ and $K_{2}$ are some constants. Setting
\begin{equation*}
\varepsilon_{k}=\rho^{k}, \ \ \rho\in (0,1),
\end{equation*}
from \eqref{oskolkovEstimate_1} it follows  the following exact order estimate 

\begin{equation*}
 {\cal E}_{n}(C(\varepsilon))_{C}=\sup\limits_{f\in C(\varepsilon) }
\| f(\cdot)-S_{n-1}(f;\cdot)\|_{C} \asymp \rho^{n},
\end{equation*}
at the same time the formuls  \eqref{SerdyukStepanyuk_FilomatLebeqIneq0} allows to write only the following estimate 
\begin{equation*}
 {\cal E}_{n}(C(\varepsilon))_{C}=o\left( \rho^{n} \log n \right).
\end{equation*}
 
 Based on the mentioned above, the problem of establishing of asymptotically unimproved Lebesgue inequalities for certain subsets of the space $C$ arose. The works \cite{Stechkin 1980, Stepanets1989N4, StepanetsSerdyuk2000} and others 
 are devoted to solving of this problem

 Investigation of asymptotic behaviour of the deviations 
 \begin{equation*}
  \rho_{n}(f; \cdot)= f(\cdot)-S_{n-1}(f;\cdot)
  \end{equation*}
  on certain functional compacts in the uniform metric was initiated in the work of Kolmogorov \cite{Kol}, who introduced the quantity 
  \begin{equation}\label{kolProblem}
 {\cal E}_{n}(\mathfrak{N})_{X}=\sup\limits_{f\in
\mathfrak{N}}\|f(\cdot)-S_{n-1}(f;\cdot)\|_{X},  \ \mathfrak{N}\subset X,
  \end{equation}
 where $\mathfrak{N}$ is some class of linear normed space  $X$.
For $X=C$ and for the class of differentiable functions  $\mathfrak{N}=W^{r}$,
  where  $W^{r}$, $r\in\mathbb{N}$ is the set of continuous functions, which have the absolutely continuous derivatives up to and including the $r$-th order, such that
$\|f^{(r)}\|_{\infty} \leq 1$,
Kolmogorov \cite{Kol}
 established the following asymptotic equality as   $n\rightarrow\infty$ 
  \begin{equation}\label{kolm1935}
 {\cal E}_{n}(W^{r}_{\infty})_{C}=\frac{4}{\pi^{2}}\frac{\ln n}{n^{r}}+\mathcal{O}\Big(\frac{1}{n^{r}}\Big),
\end{equation}
where $\mathcal{O}$ is the quantity uniformly bounded with respect to  $n$.

Later Nikolsky \cite{Nikolsky1945} established the asymptotic equality as  $n\rightarrow\infty$ for the quantities  ${\cal E}_{n}(W^{r} H^{\alpha})_{C}$,
where
  \begin{equation*}
W^{r}H^{\alpha}=\{f \in C^{r}: |f^{(r)}(x) -|f^{(r)}(x') | \leq |x-x'|^{\alpha}, \ \alpha\in(0,1), r\in \mathbb{Z}_{+]}\},
\end{equation*}
$C^{r}$ is the set of continuous functions, which  have the absolutely continuous derivatives up to and including the $r$-th order.

Namely, he showes that as  $n\rightarrow\infty$ 
 \begin{equation}\label{nikolsky_Halpha}
 {\cal E}_{n}(W^{r}H^{\alpha})_{C}= \frac{2^{\alpha+1} \ln n}{\pi^{2} n^{r+\alpha}} \int\limits_{0}^{\frac{\pi}{2}} t^{\alpha} \sin t dt+\mathcal{O}\Big(\frac{1}{n^{r+\alpha}}\Big),
\end{equation}

Nikolsky has also obtained the fundamental results regarding the research of integral norms of deviations of Fourier sums  $\| \rho_{n}(f; \cdot)\|_{L_{1}}$,  namely he has  established the analogues of the equalities \eqref{kolm1935} and \eqref{nikolsky_Halpha} in the integral metric $L_{1}$.
It should be noted that the application of Lebesgue inequality  \eqref{SerdyukStepanyuk_FilomatLebeqIneq} for the functions from the classes $W^{r}$, $r\in\mathbb{N}$ or  $W^{r}H^{\alpha}$,  $r\in\mathbb{Z}^{+}$, $\alpha\in(0,1)$ in order to obtain asymptotically accurate estimates of \eqref{kolm1935} and \eqref{nikolsky_Halpha} turns out to be ineffective.

The research of Kolmogorov and Nikolsky started a whole direction in the theory of approximations. Their results were extended to wider functional compacts $\mathfrak{N}$, to other normalized spaces $X$ and to other linear methods of polynomial approximation.
The problem about finding the asymptotic equalities as  $n\rightarrow\infty$ for the quantities
\begin{equation}\label{kolProblem_linMethod}
 {\cal E}(\mathfrak{N}, U_{n})_{X}=\sup\limits_{f\in
\mathfrak{N}}\|f(\cdot)-U_{n}(f;\cdot)\|_{X},  
  \end{equation}
where $\mathfrak{N}$ is a given class of  $2\pi$--periodic functions from the normed space  $X$, $U_{n}$ is a fixed linear polynomial method of approximation,  is called the Kolmogorov-Nikolsky problem for the class $\mathfrak{N}$ and the method  $U_{n}$ in the space $X$. At the same time, the constant in the main term of the asymptotic expansion of the quantity of the form \eqref{kolProblem_linMethod} is sometimes called as the Kolmogorov-Nikolsky constant. As it follows from above, the behavior of the Kolmogorov-Nikolsky constants for Fourier sums is not related to the behavior of the Lebesgue constants.

More details on the results about finding the solution of the Kolmogorov-Nikolsky problem for Fourier sums on classes of differential functions in the sense of Weyl-Nady $W^{r}_{\beta,p}$ and classes $(\psi,\ beta)$--differentiable functions in the sense of Stepants $L^{\psi}_{\beta,p}$ can be found in the works \cite{Kol, Korn, Nikolsky1945, Serdyuk2005Int,  Stechkin 1980, Step 1984, Stepanets1, Stepanets_Serdyuk_Shydlich, Teljakovsky1968, Teljakovsky 1989} and others.

In this work, we will limit ourselves to a review of known results for solving the Kolmogorov-Nikolsky problem for Fourier sums on the classes of generalized Poisson integrals $C^{\alpha,r}_{\beta,p}$, which are defined below.

Denote by $C^{\alpha,r}_{\beta,p}, \ \alpha>0, \ r>0, \ \beta\in
    \mathbb{R}, \ 1\leq p\leq\infty,$ the set of all  $2\pi$--periodic functions, representable for  all
$x\in\mathbb{R}$ as convolutions of the form (see, e.g., \cite[p.~133]{Stepanets1})
\begin{equation}\label{conv}
f(x)=\frac{a_{0}}{2}+\frac{1}{\pi}\int\limits_{-\pi}^{\pi}P_{\alpha,r,\beta}(x-t)\varphi(t)dt,
\ a_{0}\in\mathbb{R}, \ \varphi\in B_{p}^{0}, \
\end{equation}
$$
B_{p}^{0}=\left\{\varphi: \ ||\varphi||_{p}\leq 1, \  \varphi\perp1\right\},
 \ 1\leq p\leq \infty,
$$
with fixed generated kernels
$$
P_{\alpha,r,\beta}(t)=\sum\limits_{k=1}^{\infty}e^{-\alpha k^{r}}\cos
\big(kt-\frac{\beta\pi}{2}\big), \ \alpha>0, \ r>0, \  \beta\in
    \mathbb{R}.
$$
The kernels  $P_{\alpha,r,\beta}(t)$ are called as generalized Poisson kernels. For $r=1$ and $\beta=0$ the kernels $P_{\alpha,r,\beta}(t)$ are usual Poisson kernels of harmonic functions.

For any $r>0$ the classes  $C^{\alpha,r}_{\beta,p}$ belong to set of infinitely differentiable
 $2\pi$--periodic functions $D^{\infty}$, i.e., $C^{\alpha,r}_{\beta,p}\subset D^{\infty}$ (see, e.g., \cite[p. 128]{Stepanets1}, \cite{Stepanets_Serdyuk_Shydlich}).
For  $r\geq1$ the classes  $C^{\alpha,r}_{\beta,p}$
consist of functions  $f$,  admitting a regular extension into the strip $|\mathrm{Im} \ z|\leq c, \ c>0$ in the complex
plane (see, e.g., \cite[p.~141]{Stepanets1}), i.e., are the classes of analytic functions.
For $r>1$
the classes  $C^{\alpha,r}_{\beta,p}$  consist of functions  regular on the whole complex plane,
i.e., of entire functions (see, e.g., \cite[p.~131]{Stepanets1}). Besides,  it follows from the Theorem 1 in \cite{Stepanets_Serdyuk_Shydlich2009} that for any $r>0$ the embedding holds  $C^{\alpha,r}_{\beta,p}\subset \mathcal{J}_{1/r}$, where $\mathcal{J}_{a}, a>0,$ are known Gevrey classes
$$
\mathcal{J}_{a}=\bigg\{f\in D^{\infty}: \ \sup\limits_{k\in \mathbb{N}}\Big(\frac{\|f^{(k)}\|_{C}}{(k!)^{a}}\Big)^{1/k}<\infty \bigg\}.
$$

Approximation properties of classes of generalized Poisson integrals  $C^{\alpha,r}_{\beta,p}$ in metrics of spaces $L_{s}$, $1\leq s\leq \infty$,   were considered in  
  \cite{Serdyuk2005}--\cite{Serdyuk2012},   \cite {Serdyuk2011}--\cite {S_S2},   \cite {SerdyukStepanyuk2017},   \cite {SerdyukStepanyukAnalysis},   \cite {Step 1984},   \cite {Step monog 1987},  \cite{StepanetsSerdyuk2000}, \cite{Stepanets1},   
   \cite  {Teljakovsky 1989}--\cite{Temlyakov1990Vekya},      
 from the  viewpoint of order or asymptotic estimates for approximations by Fourier sums, best approximations and widths.

On the classes of generalized Poisson integrals $C^{\alpha,r}_{\beta,p}$,  we are interested in asymptotic equalities as  ${n\rightarrow\infty}$ of the quantities
  \begin{equation}\label{sum}
 {\cal E}_{n}(C^{\alpha,r}_{\beta,p})_{C}=\sup\limits_{f\in
C^{\alpha,r}_{\beta,p}}\|f(\cdot)-S_{n-1}(f;\cdot)\|_{C},  \ r>0, \ \alpha>0, \ \beta\in
    \mathbb{R},\ 1\leq p \leq \infty.
  \end{equation}
  


  Nikolsky \cite{Nikolsky 1946} investigated the approximation properties of the Fourier sums on some classes of analytic functions and obtained the results, from which it follows that in the case  ${r=1}$, $p=\infty$ the following asymptotic equality is true
 \begin{equation}\label{triang}
  {\cal E}_{n}(C^{\alpha,1}_{\beta,\infty})_{C}=e^{-\alpha n}\Big(\frac{8}{\pi^{2}}{\bf K}(e^{-\alpha})+\mathcal{O}(1)n^{-1}\Big),
 \end{equation}
      where
  $$
\mathbf{K}(q):=\int\limits_{0}^{\frac{\pi}{2}}\frac{dt}{\sqrt{1-q^{2}\sin^{2}t}}, \ q\in(0,1),
$$
is a complete elliptic integral of the first kind, and  $\mathcal{O}$ is a quantity uniformly bounded in parameters  $n$ and $\beta$.

  Later, the equality (\ref{triang}) was clarified by  Stechkin \cite[p. 139]{Stechkin 1980}, who established the asymptotic formula
\begin{equation}\label{stechkin}
  {\cal E}_{n}(C^{\alpha,1}_{\beta,\infty})_{C}=
 e^{-\alpha n}\Big(\frac{8}{\pi^{2}}\mathbf{K}(e^{-\alpha})+\mathcal{O}\frac{e^{-\alpha}}{(1-e^{-\alpha})n}\Big), \  \ \alpha>0, \ \beta\in\mathbb{R},
  \end{equation}
where  $\mathcal{O}(1)$   is a quantity uniformly bounded  in all analyzed parameters.

In work \cite{Serdyuk2005} for $r=1$ and arbitrary values of  $1\leq p\leq\infty$ for quantities ${\cal E}_{n}(C^{\alpha,r}_{\beta,p})_{C}$, $\alpha>0$, $\beta\in\mathbb{R}$,
the following equality was established
\begin{equation}\label{hhj}
  {\cal E}_{n}(C^{\alpha,1}_{\beta,p})_{C}=e^{-\alpha n}\bigg(\frac{2}{\pi^{1+\frac{1}{p'}}}\|\cos t\|_{p'}K(p',e^{-\alpha})+O(1)\frac{e^{-\alpha}}{n(1-e^{-\alpha})^{s(p)}}\bigg), \ \ \end{equation}
where $p'=\frac{p}{p-1}$,
$$
s(p):={\left\{\begin{array}{cc}
1, \ & p=\infty,  \\
2, \ & p\in[1,2)\cup(2,\infty), \\
-\infty, &
p=2, \
  \end{array} \right.}
$$
$$
K(p',q):=\frac{1}{2^{1+\frac{1}{p'}}}\Big\|(1-2q\cos t+q^{2})^{-\frac{1}{2}} \Big\|_{p'}, \ q\in(0,1),
$$
and $\mathcal{O}(1)$ is a quantity uniformly bounded in $n$, $p$, $\alpha$ and $\beta$.
For $p=\infty$, by virtue of  the known equality  $K(1,q)=\mathbf{K}(q)$ (see, e.g., formula 3.674(1) of \cite[p. 401]{Gradshteyn}), the estimate (\ref{hhj}) coincides with the estimate   (\ref{stechkin}).

In \cite{Serdyuk2012} it was proved that in the case $1\leq p'<\infty$ the following equality takes place 
\begin{equation}\label{elliptic_hypergeometric}
K(p',q)=\frac{\pi^{\frac{1}{p'}}}{2}F^{\frac{1}{p'}}\Big(\frac{p'}{2}, \frac{p'}{2}; 1; q^{2}\Big),  \ q\in(0,1),
\end{equation}
where $F(a,b;c;d)$ is Gauss hypergeometric function
  $$
F(a,b;c;z)=1+\sum\limits_{k=1}^{\infty}\frac{(a)_{k}(b)_{k}}{(c)_{k}}\frac{z^{k}}{k!},
$$
$$
(x)_{k}:=x(x+1)(x+2) ...(x+k-1).
$$
   
So, taking into account the equality     \eqref{elliptic_hypergeometric}, where $q=e^{-\alpha}$, the estimate \eqref{hhj} can be written in the form 
   \begin{equation}\label{hhj_v1}
  {\cal E}_{n}(C^{\alpha,1}_{\beta,p})_{C}=e^{-\alpha n}\bigg(\frac{\|\cos t\|_{p'} }{\pi}F^\frac{1}{p'}\left(\frac{p'}{2}, (\frac{p'}{2} ; 1; e^{-2\alpha} \right)+O(1)\frac{e^{-\alpha}}{n(1-e^{-\alpha})^{s(p)}}\bigg).
  \end{equation}

Note that for $p=2$ and $r=1$ formula  (\ref{hhj}), and also \eqref{hhj_v1} becomes the equality
\begin{equation}\label{serd2005_eq}
 {\cal E}_{n}(C^{\alpha,1}_{\beta,2})_{C}=\frac{1}{\sqrt{\pi(1-e^{-2\alpha})}}e^{-\alpha n}, \ \alpha>0, \ \beta\in\mathbb{R}, \ n\in\mathbb{N},
 \end{equation}
(see \cite{Serdyuk2005}).
Moreover,  it follows from  \cite{Serdyuk2011} that for $p=2$  for the quantities
${\cal E}_{n}(C^{\alpha, r}_{\beta,p})_{C}$ the following equalities take place
\begin{equation}\label{serd2011}
 {\cal E}_{n}(C^{\alpha,r}_{\beta,2})_{C}=\frac{1}{\sqrt{\pi}}\Big(\sum\limits_{k=n}^{\infty}e^{-2\alpha k^{r}}\Big)^{\frac{1}{2}}, \ \alpha>0,  \ r>0,  \ \beta\in\mathbb{R}, \ n\in\mathbb{N}.
\end{equation}

It is clear that for $r=1$ the formula \eqref{serd2011} transforms into \eqref{serd2005_eq}.

 In the case of ${r>1}$ and   $p=\infty$ the asymptotic equalities for the quantities  ${\cal E}_{n}(C^{\alpha,r}_{\beta,p})_{C}$, $\alpha>0$, $\beta\in\mathbb{R}$,   were obtained by Stepanets
 \cite[Chapter 3, Section~9]{Step monog 1987}, who showed that for any  $n\in\mathbb{N}$
\begin{equation}\label{step_prepr}
  {\cal E}_{n}(C^{\alpha,r}_{\beta,\infty})_{C}=\Big(\frac{4}{\pi}+\gamma_{n}\Big)e^{-\alpha n^{r}},
\end{equation}
where
$$
|\gamma_{n}|<2\Big(1+\frac{1}{\alpha r n^{r-1}}\Big)e^{-\alpha r n^{r-1}}.
$$

 Later Telyakovsky  \cite{Teljakovsky 1989} clarified  the asymptotic equality \eqref{step_prepr} by improving the estimate for the residual term. Namely, it follows from \cite{Teljakovsky 1989} that
\begin{equation}\label{tel}
  {\cal E}_{n}(C^{\alpha,r}_{\beta,\infty})_{C}=\frac{4}{\pi}e^{-\alpha n^{r}}+
  \mathcal{O}(1)\Big(e^{-\alpha ( 2(n+1)^{r}-n^{r})}
  +\Big(1+\frac{1}{\alpha r(n+2)^{r}}\Big)e^{-\alpha (n+2)^{r}}  \Big),
\end{equation}
where  $\mathcal{O}(1)$ is a quantity uniformly bounded in all analyzed parameters.

For $r>1$ and for arbitrary values of  $1\leq p\leq\infty$ the asymptotic equalities for the quantities ${\cal E}_{n}(C^{\alpha,r}_{\beta,p})_{C}$, $\alpha>0$, $\beta\in\mathbb{R}$,
are found in \cite{Serdyuk2005} and have the form
\begin{equation}\label{ser}
  {\cal E}_{n}(C^{\alpha,r}_{\beta,p})_{C}=e^{-\alpha n^{r}}\Big(\frac{\|\cos t\|_{p'}}{\pi}+\mathcal{O}(1)\Big(1+\frac{1}{\alpha r n^{r-1}}\Big)e^{-\alpha n^{r-1}}\Big),
\end{equation}
where $\mathcal{O}$ is a quantity uniformly bounded in all analyzed parameters. For ${p=\infty}$ the formula (\ref{ser}) follows from (\ref{step_prepr}) and
(\ref{tel}).

 Concerning the case ${0<r<1}$,  except the presented above case  $p=2$, asymptotic equalities for quantities  ${\cal E}_{n}(C^{\alpha,r}_{\beta,p})_{C}$,  $\alpha>0$,  $\beta\in\mathbb{R}$, were known only for $p=\infty$ due to the work of  Stepanets  \cite{Step 1984}, who showed that
\begin{equation}\label{step}
  {\cal E}_{n}(C^{\alpha,r}_{\beta,\infty})_{C}=\frac{4}{\pi^{2}}e^{-\alpha n^{r}}\ln n^{1-r}+\mathcal{O}(1)e^{-\alpha n^{r}}, \ r\in(0,1), \ \beta\in\mathbb{R},
\end{equation}
where $\mathcal{O}(1)$ is a quantity uniformly bounded in $n$ and $\beta$.

     In case of  ${0<r<1}$ and  $1\leq p<\infty$ the following order estimates for quantities ${\cal E}_{n}(C^{\alpha,r}_{\beta,p})_{C}$,  $\alpha>0$,  $\beta\in\mathbb{R}$, hold  (see, e.g., \cite{Temlyakov1990Vekya}, \cite{S_S})
\begin{equation}\label{order}
{\cal E}_{n}(C^{\alpha,r}_{\beta,p})_{C}\asymp e^{-\alpha n^{r}}n^{\frac{1-r}{p}}.
\end{equation}

We remark that  for ${0<r<1}$ and $1\leq p<\infty$  Fourier sums provide the order of the best approximations of classes $C^{\alpha,r}_{\beta,p}$,  $\alpha>0$,  $\beta\in\mathbb{R}$, in the uniform metric, i.e. (see, e.g., \cite{S_S}, \cite{S_S2})
$$
{\cal E}_{n}(C^{\alpha,r}_{\beta,p})_{C}\asymp{ E}_{n}(C^{\alpha,r}_{\beta,p})_{C}\asymp e^{-\alpha n^{r}}n^{\frac{1-r}{p}},
$$
where
$$
{ E}_{n}(C^{\alpha,r}_{\beta,p})_{C}=\sup\limits_{f\in
C^{\alpha,r}_{\beta,p}}\, \inf\limits_{t_{n-1}\in\mathcal{T}_{2n-1}}\|f-t_{n-1}\|_{C},
$$
and $\mathcal{T}_{2n-1}$ is the subspace of all trigonometric polynomials $t_{n-1}$ of degree not higher than ${n-1}$.

On the one hand, this fact encourages to research more deeply approximation properties of Fourier sums in given situations, and on the other hand it separates the case $1\leq p<\infty$ from  considered earlier case $p=\infty$, where order equality ${\cal E}_{n}(C^{\alpha,r}_{\beta,p})_{C}\asymp{E}_{n}(C^{\alpha,r}_{\beta,p})_{C}$, $0<r<1$, doesn't take place.

Besides, as  follows from Temlyakov work   \cite{Temlyakov1990Vekya} for $2\leq p<\infty$ the quantities of approximations by Fourier sums provide the order of the linear widths  $\lambda_{2n}$ (definition of  $\lambda_{m}$ see, e.g., \cite[Chapter 1, Section 1.2]{Korn}) of  the classes $C^{\alpha,r}_{0,p}$, i.e.
$$
\lambda_{2n}(C^{\alpha,r}_{0,p},C)\asymp{\cal E}_{n}(C^{\alpha,r}_{0,p})_{C}.
$$

In the paper \cite{SerdyukStepanyukAnalysis} the authors established the asymptotically sharp estimates of the quantities ${\cal E}_{n}(C^{\alpha,r}_{\beta,p})_{C}$,  ${\alpha>0}$,  $\beta\in\mathbb{R}$,
 for any $0< r<1$ and $1\leq p\leq\infty$. In particular, it is proved, that for  $r\in(0,1)$, $\alpha>0$, $\beta\in\mathbb{R}$ and $1\leq p\leq\infty$ as $n\rightarrow\infty$ the following asymptotic equalities take place
 $$
 {\cal E}_{n}(C^{\alpha,r}_{\beta,p})_{C}=e^{-\alpha n^{r}}n^{\frac{1-r}{p}}\bigg(
\frac{\|\cos t\|_{p'}}{\pi^{1+\frac{1}{p'}}(\alpha r)^{\frac{1}{p}}}F^{\frac{1}{p'}}\Big(\frac{1}{2}, \frac{3-p'}{2}; \frac{3}{2}; 1\Big)+
$$
\begin{equation}\label{tan}
+\frac{\mathcal{O}(1)}{n^{\min\{\frac{1-r}{p},r\}} }\bigg), \ \ 1< p\leq\infty, \ \ \frac{1}{p}+\frac{1}{p'}=1,
  \end{equation}
  \begin{equation}\label{tan1}
{\cal E}_{n}(C^{\alpha,r}_{\beta,1})_{C}=
e^{-\alpha n^{r}}n^{1-r}\bigg(\frac{1}{\pi\alpha r}+\frac{\mathcal{O}(1)}{n^{\min\{1-r,\  r\}} }\bigg),
\end{equation}
   where  $\mathcal{O}(1)$ is a quantity uniformly bounded with respect to $n$ and $\beta$.
   
Formulas (\ref{tan}) and (\ref{tan1}) together with formulas (\ref{triang})--(\ref{hhj}), (\ref{step_prepr})--(\ref{step}) give the solution of Kolmogorov--Nikolsky problem about strong asymptotic of quantities (\ref{sum}) as $n\rightarrow\infty$ for all admissible values of parameters of problem.   
   
 Herewith, in the paper \cite{SerdyukStepanyukAnalysis} the authors found the estimates for remainders in (\ref{tan}) and (\ref{tan1}),  which are expressed via absolute values and the parameters of the problem $\alpha, r, p$ in the explicit form. These estimates can be used for practical application, since they allow effectively to estimate the errors of uniform approximations of functions from the classes $C^{\alpha,r}_{\beta,p}$ by their partial Fourier sums.

The following table  contains the exact values  of constants (Kolmogorov--Nikolsky constants) in the main term $A_n$ of asymptotic expansion of quantities $ {\cal E}_{n}(C^{\alpha,r}_{\beta,p})_{C}$ of the form
 $$
  {\cal E}_{n}(C^{\alpha,r}_{\beta,p})_{C}=e^{-\alpha n^{r}}(A_n+o(A_n)).
 $$

\begin{center}\label{Tabl}
\begin{tabular}{|c|c|c|c|c|}
	\hline
	\multicolumn{2}{|c|}{}  & \multicolumn{3}{|c|}{ $r$}  \\
\cline{3-5}
\multicolumn{2}{|c|}{\raisebox{0.2ex}[0cm][0cm]{{ $A_n$}}  } &  \raisebox{-0.9ex}[0cm][0cm]{\boldmath{ $(0,1)$} }&  \raisebox{-0.9ex}[0cm][0cm]{\boldmath{$1$}} & \raisebox{-0.9ex}[0cm][0cm]{ \boldmath{$(1,\infty)$}} \\
\multicolumn{2}{|c|}{ } &   &  &   \\
\hline
\ &  &    &   &  \\
 \ &  \raisebox{-2.8ex}[0cm][0cm]{  \boldmath{$\infty$}  } & \cite{Step 1984} Stepanets (1984)     & \cite{Nikolsky 1946} Nikolsky (1946)  & \cite{Step monog 1987} Stepanets   \\
 \  &  &      &   & (1987)  \\ 
                    &   &   & \cite{Stechkin 1980} Stechkin (1980)  & \cite{Teljakovsky 1989}  Telyakovskii  \\
                     &  & $\frac{4}{\pi^{2}}(1-r)\ln n$ & \ &  (1989) \\
                      &   &                                & $\frac{8}{\pi^{2}}\mathbf{K}(e^{-\alpha})$ &  $\frac{4}{\pi}$ \\
\ &  &    &   &  \\
\cline{2-5}
\ &  &    &   &  \\
 \raisebox{-1.4ex}[0cm][0cm]{  \boldmath{$p$}  }&  \raisebox{-3.8ex}[0cm][0cm]{  \!\!\! \!\boldmath{$(1,  \infty)$} \!\!\!\!} &  \raisebox{0.8ex}[0cm][0cm]{ \cite{SerdyukStepanyukAnalysis} Serdyuk Stepanyuk (2019)}  &  \raisebox{0.8ex}[0cm][0cm]{ \cite{Serdyuk2005}  Serdyuk (2005) }  & \raisebox{0.8ex}[0cm][0cm]{ \cite{Serdyuk2005}  Serdyuk  } \\
   &  \  &    &   &  (2005)   \\ 
                        &  &   $\frac{\|\cos t\|_{p'} F^{\frac{1}{p'}}(\frac{1}{2}, \frac{3-p'}{2}; \frac{3}{2}; 1)}{\pi^{1+\frac{1}{p'}}(\alpha r)^{\frac{1}{p}}} n^{\tfrac{1-r}{p}}$  & $\frac{\|\cos t\|_{p'}}{\pi} F^{\frac{1}{p'}}(\frac{p'}{2}, \frac{p'}{2}; 1; e^{-2\alpha})$
                            & $\frac{\|\cos t\|_{p'}}{\pi}$ \\
                            \ &  &    &   &  \\
 \cline{2-5}
 &   &    &     &    \\
   &  \raisebox{-2.8ex}[0cm][0cm]{  \boldmath{$1$} } &   \cite{SerdyukStepanyukAnalysis} Serdyuk Stepanyuk (2019) & \cite{Serdyuk2005} Serdyuk (2005)    & \cite{Serdyuk2005} Serdyuk   \\
&  &     &     &   (2005)   \\   
                        &  &   \raisebox{-0.9ex}[0cm][0cm]{$\frac{1}{\pi\alpha r} n^{1-r}$ } & \raisebox{-0.9ex}[0cm][0cm]{$\frac{1}{\pi(1-e^{-\alpha})} $} & \raisebox{-0.9ex}[0cm][0cm]{$\frac{1}{\pi}$ }\\
  \tiny{ \ }  &  \tiny{ \ }  &   \tiny{ \ } & \tiny{ \ } &  \tiny{ \ } \\
 \hline
\end{tabular}
\end{center}


Let us make a few remarks about the solution of the Kolmogorov-Nikolsky problem in the metrics of spaces $L_{p}$ .

Currently the asymptotic equalities as $n\rightarrow\infty $ of the quantities
\begin{equation}\label{sumIntegr}
 {\cal E}_{n}(C^{\alpha,r}_{\beta,1})_{L_{s}}=\sup\limits_{f\in
C^{\alpha,r}_{\beta,1}}\|f(\cdot)-S_{n-1}(f;\cdot)\|_{s},  \ r>0, \ \alpha>0, \ \beta\in
    \mathbb{R},\ 1\leq s \leq \infty,
  \end{equation}
are known for all admissible values ??of the problem parameters. At the same time, it turned out that the following asymptotic equality takes place
\begin{equation}\label{asympIntegr}
 \lim\limits_{n\rightarrow \infty}\frac{ {\cal E}_{n}(C^{\alpha,r}_{\beta,p})_{C}}{ {\cal E}_{n}(C^{\alpha,r}_{\beta,1})_{L_{p'}}} =1, \ \ \frac{1}{p}+ \frac{1}{p'} =1,
  \end{equation}
  which holds for all $\alpha>0$, $r>0$ and  $\beta\in \mathbb{R}$, $1\leq p \leq \infty$. The works \cite{Serdyuk2005Int, S_S2, Step 1984, Step monog 1987, Stepanets1} are devoted to the establishment of asymptotic equalities for quantities of the form \eqref{sumIntegr}.
  
  In the works \cite{MusienkoSerdyuk_UMJ1, MusienkoSerdyuk_UMJ2, MusienkoSerdyuk_Zb_IM2010, SerdyukStepanyukFilomat, SerdyukStepanyuk_Jaen, Stepanets1989N4,StepanetsSerdyuk2000} for the functions $f$ from the sets $C^{\alpha,r}_{\beta}L_{p}$ or $C^{ \alpha,r}_{\beta}C$, defined by convolutions of the form \eqref{conv},
 in which the functions $\varphi$ belong to the spaces $L_{p}$ or $C$, respectively, the  asymptotically best possible Lebesgue-type inequalities were established. In these equalities the uniform norms of the deviations of Fourier sums $\| f-S-{n-1}\|_{C}$ were estimated via the best approximations of the functions $\varphi$ in the metrics of spaces $L_{p}$ or $C$.





\end{document}